\documentstyle[12pt]{article}

\def\rr{{\bf R}}
\def\qq{{\bf Q}}
\def\res{\upharpoonright}
\def\qed{\nopagebreak\par\noindent\nopagebreak$\blacksquare$\par}
\input amssym.def
\input amssym

\def\al{\alpha}
\def\be{\beta}
\def\ga{\gamma}
\def\de{\delta}
\def\ep{\epsilon}
\def\la{\lambda}
\def\om{\omega}

\def\su{\subseteq}
\def\se{\setminus}
\def\pf{{\par\bf Proof.\ \ }}
\def\cantor{2^\omega}
\def\cross{\times}
\def\meas{{\Bbb B}}

\newtheorem{theorem}{Theorem}
\newtheorem{lemma}[theorem]{Lemma}

\begin{document}

\begin{center}
{\bf\Large
Measurability of functions with approximately
continuous vertical sections and measurable horizontal sections}
\end{center}

\bigskip
\begin{center}
  By M. Laczkovich and Arnold W. Miller
\end{center}
\bigskip

A function $f:\rr\to\rr$ is approximately continuous iff
it is continuous in the density topology, i.e., for any
open set $U\su\rr$ the set $E=f^{-1}(U)$ is measurable
and has Lebesgue density one at each of its points.
Approximate continuity was introduced by Denjoy \cite{den} in
his study of derivatives.  Denjoy proved that
bounded approximately continuous functions are
derivatives.  It follows from this that approximately continuous
functions are Baire 1, i.e., pointwise limits of continuous
functions.  For more on these concepts, see Bruckner \cite{bru},
Luke\v{s}, Mal\'{y}, Zaji\v{c}ek \cite{lmz}, Tall \cite{tall},
and Goffman, Neugebauer, Nishiura \cite{gnn}.

For any $f:\rr^2\to\rr$ define
$$f_x(y)=f^y(x)=f(x,y)$$
for any $x,y\in\rr$.  A function $f:\rr^2\to\rr$ is separately
continuous if $f_x$ and $f^y$ are continuous for every $x,y\in\rr$.
Lebesgue \cite{leb} in his first paper proved that any
separately continuous function is Baire 1. He also
showed  that if $f_x$ is continuous for all $x$ and $f^y$
Baire $\al$ for all $y$, then $f$ is Baire $\al+1$ (see
Kuratowski \cite{kur} p. 378).  For more historical comments and
 generalizations see Rudin \cite{rud}.
Sierpi\'{n}ski
\cite{sier} showed that there exists a nonmeasurable $f:\rr^2\to\rr$
which is separately Baire 1. (The characteristic function of
a nonmeasurable subset of the plane which meets every horizontal
and vertical line in at most one point.)

In this paper we shall prove:

\begin{theorem}\label{thm1}
Let $f:\rr^2 \to \rr$ be such that
$f_x$ is approximately
continuous and $f^y$ is Baire 1 for every $x,y\in \rr .$
Then $f$ is Baire 2.
\end{theorem}

\begin{theorem}\label{thm1.5}
Suppose there exists a real-valued measurable cardinal.
Then for any function $f:\rr^2\to\rr$ and $\al<\om_1$,
if $f_x$ is approximately continuous
and $f^y$ is Baire $\al$ for every $x,y\in\rr$, then
$f$ is Baire $\alpha+1$ as a function of two variables.
\end{theorem}

\begin{theorem}\label{thm2}
{\rm (i)} Suppose that $\rr$ can be covered by $\om _1$
closed null sets. Then there exists a nonmeasurable
function $f:\rr^2 \to \rr$ such that $f_x$ is approximately
continuous and $f^y$ is Baire 2 for every $x,y\in \rr .$
\par\noindent {\rm (ii)}
Suppose that $\rr$ can be covered by $\om _1$
null sets. Then there exists a nonmeasurable
function $f:\rr^2 \to \rr$ such that $f_x$ is approximately
continuous and $f^y$ is Baire 3 for every $x,y\in \rr .$
\end{theorem}

\begin{theorem}\label{thm3}
In the random real model for any function $f:\rr^2\to\rr$
if $f_x$ is approximately continuous
and $f^y$ is measurable for every $x,y\in\rr$, then
$f$ is measurable as a function of two variables.
\end{theorem}

{\bf Remarks.} Davies \cite{dav} showed
that any function of two variables which is separately
approximately continuous is Baire 2.
Theorem \ref{thm1} which generalizes this was announced in
Laczkovich and Petruska \cite{4},
but the proof was never published.
In  Davies and  Draveck\'y \cite{2} and Grande \cite{3} it
is shown that CH implies the existence of
a nonmeasurable function $f$ such that $f_x$ is approximately
continuous for every $x$ and $f^y$ is measurable for every $y$.
It was pointed out in Laczkovich and Petruska \cite{4} on p. 170
that the construction, in fact,
gives Baire 2 sections. Our Theorem \ref{thm2} is a refinement
of this observation.  Note that Bartoszynski and Shelah \cite{barshe}
have shown that it is relatively consistent with ZFC
that $\rr$ is the union of $\omega_1$ meager null sets, but
not the union of $\omega_1$ closed null sets.  It is well
known that $\rr$ can be the union of $\omega_1$ closed null sets
and the continuum arbitrarily large.

In Theorem \ref{thm1.5} we only use that for any family of
continuum many subsets of the real line there exists
a measure extending Lebesgue measure and making the family
measurable.  This is slightly weaker than a real-valued
measurable and has the consistency strength of a weakly
compact cardinal (see Carlson \cite{carl}).

It follows from
Lebesgue's argument that any function $f:\rr^2\to \rr$ such
that $f_x$ is continuous and $f^y$ is measurable for all $x,y\in\rr$
must be measurable as a function of two variables.
Theorems \ref{thm2} and \ref{thm3} show that this fact
is independent of set theory if we replace continuous by approximately
continuous.

\bigskip
\bf Proof of Theorem \ref{thm1}. \rm
This is an immediate
consequence of the following theorem due to Bourgain, Fremlin
and Talagrand \cite{1}.
\begin{theorem}[Bourgain, Fremlin, Talagrand] \label{BFT}
Let $(X,\Sigma , \mu )$
be a probability space
and let $f:X\times \rr \to \rr$ be bounded. If $f_x$ is Baire 1 for
every $x\in X$ and $f^y$ is measurable for every $y\in\rr$,
then the function $$y\mapsto \int_X f^y d\mu(x) \ (y\in \rr )$$
is Baire 1.
\end{theorem}

Suppose that $f_x$ is approximately continuous and $f^y$ is Baire 1 for
every $x,y\in \rr .$  Without loss of generality we may assume
that $f$ is bounded. (Otherwise, let $h:\rr\to (0,1)$ be
a homeomorphism.  Then $h \circ f$ is approximately continuous
when $x$ is fixed and measurable when $y$ is fixed. Hence
$h\circ f$ is Baire 2 and therefore $h^{-1}\circ h \circ f=f$ is
Baire 2.)

It follows from  Theorem \ref{BFT}, that for every fixed $y,$
the function
$$x\mapsto \int_0^y f_x dt \ (x\in \rr )$$ is Baire 1.

This implies that the function
$$F(x,y)=\int_0^y f_x dt$$
is Baire 1,
since $F^y$ is Baire 1 and the family $\{ F_x : x\in \rr \}$
is uniformly continuous (in fact, uniformly Lipschitz).
The proof is this. Let $F_n :\rr^2 \to \rr$ be the function such that
$F_n(x,i/n)=F(x,i/n)$ for every $x\in \rr$ and every integer $i,$
and let $F_n (x_0 ,y)$ be linear in $y\in [(i-1)/n ,i/n]$
for every integer $i$ and every fixed $x_0.$  Then $F_n$ is Baire 1.
Indeed, let $F(x,i/n)=\lim_{j\to \infty} g_{i,j}(x),$
where $g_{i,j}:\rr\to \rr$ continuous. Let $G_j (x,i/n)=g_{i,j} (x),$
let $G_j$ be continuous in $y$ and linear for $y\in [(i-1)/n,i/n]$
for every fixed $x.$ Then $G_j$ is continuous and $G_j \to F_n ,$
so that $F_n$ is Baire 1. Finally, $F_n \to F$ uniformly,
so that $F$ is Baire 1 (see Kuratowski \cite{kur} p. 386).

Finally, since
$$f(x,y)=\lim_{n\to \infty} {{F(x,y+(1/n))-F(x,y)}\over{ 1/n}},$$
it follows that $f$ is Baire 2.
\qed
\bigskip

\bf Proof of Theorem \ref{thm1.5}. \rm
This is the same as the proof of Theorem \ref{thm1} except
we use the following generalization of the Bourgain-Fremlin-Talagrand
Theorem \ref{BFT}:
\begin{lemma}
Let $(X,\Sigma , \mu )$
be a probability space such that every subset of $X$ is in $\Sigma$
and let $f:X\times \rr \to \rr$ be bounded. For
$\al<\om_1$ if $f_x$ is Baire $\al$ for
every $x\in X$,
then the function $$F(y)=\int_X f^y d\mu(x) \mbox{ for }y\in \rr $$
is Baire $\al$.
\end{lemma}
\pf
This is proved by induction on $\al$. If $\alpha =0; $
that is, if $f_x$ is continuous for every $x,$ then the
continuity of $F$ follows from the dominated convergence theorem.
For $\al>0$,
let $\be_n$ be
a nondecreasing sequence of ordinals such that
$\sup_{n\in\om}(\be_n+1)=\al$.  Let
$\langle f_n:n\in\om \rangle$ be a sequence of uniformly bounded
functions such
that $(f_n)_x$ is Baire $\be_n$ for each $n$ and
$$\lim_{n\to\infty}f_n(x,y)=f(x,y).$$
Then by induction the function
$$F_n(y)=\int_X f_n^y d\mu(x)$$
is Baire $\be_n$.  By the dominated convergence theorem
$$\lim_{n\to\infty}F_n(y)=F(y)$$
is Baire $\al$.
\qed
Since there is a real-valued measurable cardinal we
can find an extension $\mu$ of Lebesgue measure $\la$
which makes every set of reals measurable.  The
rest of the proof is the same as Theorem \ref{thm1}.
\qed

\bigskip
\bf Proof of Theorem \ref{thm2}. \rm
Let $\rr =\cup_{\al <\om _1} C_\al ,$ where $C_\al$ is
a closed set of measure zero
for every $\al <\om _1 .$

By a Lemma of Zahorski \cite{zah} (see also
Bruckner \cite{bru} p. 28) for any $G_\de$ measure zero set $G\su\rr$
there
exists an approximately continuous $g:\rr\to [0,1]$ such that
$g^{-1}\{0\}=G$.
So for each $\al$ let $g_\al :\rr \to [0,1]$
be an approximately continuous function such that
$g_\al^{-1}\{0\}$ is a measure zero set covering
$\bigcup _{\be<\al} C_\be$.
We define
$f(x,y)=g_\al (y),$ where
$\al$ is the smallest ordinal such that $x\in C_\al .$

Obviously, $f_x$ is approximately continuous for every $x.$
For any fixed $y,$ let $\al$ be such that $y\in C_\al .$
If $x\notin \bigcup_{\be<\al}C_\al$,
then $f(x,y)=0.$  It is also clear that $f^y$ is constant on each
of the $G_\de$ sets $C_\be \se \cup_{\ga <\be} C_\ga $.  It follows
that $f^y$ is Baire 2, since the range of $f^y$ is countable
and the preimage of any set is a countable union of
$G_\de$-sets.
Finally, $f$ is not measurable, since
$$\int_\rr \biggl ( \int_\rr f_x dy \biggr ) dx >0=
\int_\rr \biggl ( \int_\rr f^y dx \biggr ) dy.$$

\bigskip
For the second part, let ${\bf R} =\cup_{\alpha <\omega _1} C_\alpha ,$ where
$\lambda (C_\alpha )=0$ for every $\alpha <\omega _1 .$ We may assume that
each $C_\alpha$ is a $G_\delta$ set. Following the proof of (i), we obtain
a nonmeasurable function $f$ such that $f_x$ is approximately continuous
for every $x.$ Also, for every $y,$ the preimage of any set by $f^y$ is a
countable union of $F_{\sigma \delta}$ sets, and thus $f^y$ is Baire 3.
\qed

\bf Proof of Theorem \ref{thm3}. \rm
We will use the following lemmas.
For a set in the plane $H\su\rr\cross\rr$ and $x,y\in\rr$ let
$$H_x=\{y\in\rr : (x,y)\in H\} \mbox{ and }
H^y=\{x\in\rr : (x,y)\in H\}.$$

\begin{lemma}\label{equiv}
The following statements are equivalent.
\begin{description}
 \item[{\rm (i)}] There exists a nonmeasurable
     function $f:\rr^2 \to \rr$ such that $f_x$ is approximately
     continuous and $f^y$ is measurable for every $x,y\in \rr$.
 \item[{\rm (ii)}] There exists a set $H\su \rr ^2$ such that
     $\la (H^y )=0$ for every $y\in\rr$, but the set
     $\{ x: \la (\rr \se H_x )=0 \}$ has positive outer measure.
\end{description}
\end{lemma}

\pf
(ii)$\Longrightarrow$(i): Suppose (ii) and let
$A=\{ x: \la (\rr \se H_x )=0 \}$. For every $x\in A$
there is a $G_\de$ null set $B_x \su \rr$ such that $\rr \se H_x \su B_x .$
This implies by Zahorski's Lemma
that for every $x\in A$ there exists an approximately
continuous function $g_x :\rr \to \rr$ such that $g_x (y)=0$ if
$y\in B_x$ and $0<g_x (y) \le 1$ if $y\notin B_x .$

For every $y\in \rr$
we define $f(x,y)=g_x (y)$ if $x\in A$, and $f(x,y)=0$ if $x\notin A.$
Then $f_x$ is approximately continuous for every $x$.  Also, $f^y$
is measurable for every $y,$ since $f^y  (x)=0$ for a.e. $x.$ Indeed,
$$f^y (x)\ne 0 \Longrightarrow x\in A,\ y\notin B_x ,
\Longrightarrow y\in H_x \Longrightarrow x\in H^y$$
and hence
$$\la (\{x: f^y (x)\ne 0 \} )\le \la (H^y )=0.$$

This implies that
$$\int_\rr \left(\;\int_\rr f^y dx \right)dy =0.$$
On the other hand,
$$\int_\rr \left(\;\int_\rr f_x dy \right)dx >0,$$
since $\int_\rr f_x dy >0$
for every $x\in A$ and $A$ has positive outer measure.
Therefore $f$ cannot be measurable.

\medskip
(i)$\Longrightarrow$(ii): Suppose (i); we may also assume that $f$
is bounded.

Since every approximately continuous function is Baire 1,
it follows as in the proof of
Theorem \ref{thm1}, that the function
$$F(x,y)=\int_0^x f^y dt$$
is Baire 1. Let
$$g(x,y)=\left\{\begin{array}{ll}
    \lim_{n\to \infty} n\cdot \bigl ( F(x+(1/n), y) -F(x,y) \bigr )&
    \mbox{if this limit exists}\\
     0 & \mbox{if it does not.} \\
                \end{array}\right.
$$
Then $g$ is Borel measurable, and for every
fixed $y,$ we have $g(x,y)=f(x,y)$ for a.e. $x$ by Lebesgue's
classical theorem.

\medskip
{\bf Claim.}  For any $g:\rr^2\to\rr$ measurable, there exists a Borel
set $B\su \rr^2$ such that $\la _2 (B)=0$ and for every
$(x,y)\notin B$
the function $g_x$ is approximately continuous at $y$.
\pf
This easily follows from the fact that if $E\su \rr ^2$ is measurable
then there is a Borel
set $B\su \rr ^2$ such that $\la _2 (B)=0$ and
$y$ is a density point of $E_x$ for every
$(x,y)\in E\se B$; see the argument on pp. 130-131
of Saks \cite{5}.  For the convenience of the reader we sketch
the proof here.  Without loss of generality, we may assume $E$ is
compact.  Fix $\ep>0$ and define
$$A_n^\ep=\{(x,y)\in E: \la(E_x\cap I)\geq(1-\ep)\la(I)
\mbox{ whenever $y\in I$ and $|I|<1/n$}\}.$$

(We use $I$ to range over nondegenerate closed intervals.)
Then it can be shown that $A_n^\ep$ is closed since $E$ is.
Therefore
$$N_\ep=E\se \bigcup_{n\in\om}A_n^\ep$$
is measurable.  By the Lebesgue density theorem, $(N_\ep)_x$ has
measure zero for every $x$ and hence by Fubini's Theorem
$N_\ep$ has planar measure zero.  Let
$$B=\bigcup_{\ep>0}N_\ep.$$
Then $\la_2(B)=0$ and $y$ is a
density point of $E_x$ for every $(x,y)\in E\se B$.
To obtain the result for $g$ let $B$ be a measure zero
subset of the plane such that for every $U$ in some countable basis
for $\rr$ if $(x,y)\in g^{-1}(U)\se B$, then
$y$ is a
density point of $(g^{-1}(U))_x=g_x^{-1}(U)$.  It follows
that $g_x$ is
approximately continuous at $y$
for every $(x,y)\in \rr\se B$.  This proves the Claim.
\qed

Let
       $$K= \{ (x,y): g(x,y)\ne f(x,y) \} ;$$
then $\la (K^y )=0$
for every $y.$ Let $x$ be fixed. Then, for $y\notin B_x ,$ the
functions $f_x$ and $g_x$ are both approximately continuous at $y$.
Therefore, if $(x,y)\in K$ then the set
       $$K_x =\{ y: f_x (y)\ne g_x (y)\}$$
is measurable and of positive measure. (This is because if two
functions are approximately continuous at a point $x$ and take
on different values there,
then there exists a measurable set with density one at $x$ where
they differ.)

Hence for any $x$, $K_x$ is measurable, and
either $K_x \su B_x$ or $K_x$ has positive measure.
Let $A=\{ x: \la (K_x )>0\} ,$ then $K\su B\cup (A\times \rr ).$ If $\la (A)
=0$ then $\la _2 (K)=0$ and $f=g$ almost everywhere, contradicting
our assumption that $f$ is not measurable. Thus $A$ has positive outer
measure.

Now, putting $H=\{ (x,y+r): (x,y)\in K,\; r\in {\qq} \} ,$ we obtain a set
such that $\la (H^y )=0$ for every $y$ and $\la (\rr \se H_x )=0$
for  $x\in A;$ and hence (ii) holds.
\qed

By the random real model we refer to any model of set theory
which is a generic extension of a countable transitive ground model
of CH by adding $\omega_2$ random reals, i.e., forcing
with the measure algebra on $2^{\om_2}$.

\begin{lemma}\label{randmod}
In the random real model the following two facts hold:
\begin{enumerate}
  \item $\rr$ is not the union of $\om_1$ measure zero
  sets.
  \item Any $Y\su\rr$ with positive outer measure
  contains a subset $Z\su Y$ of cardinality $\om_1$
  with positive outer measure.
\end{enumerate}
\end{lemma}

\pf
Lemma \ref{randmod}.1 is due to Solovay \cite{S} and is also proved
in Kunen \cite{K} 3.18 and probably Jech \cite{J}.
Lemma \ref{randmod}.2 is probably due to Kunen
(see remark in Tall \cite{tall} p. 283),
but we don't know of a published proof, so we include one here.

Since $\cantor$ and $[0,1]$ are measure isomorphic, we may work in
$\cantor$.
For any set $\Sigma$ let $2^\Sigma$ be product space of the two point set
$2=\{0,1\}$ with the usual product measure and topology.
Let $\meas(\Sigma)$ denote
the measure algebra, i.e., the Borel subsets of $2^\Sigma$ modulo the
measure zero sets.  This is a complete
boolean algebra which satisfies the countable chain condition.

Let $M$ be a countable standard model of ZFC+CH.
For any
set $\Sigma$ in $M$ let $\meas(\Sigma)^M$ denote
the measure algebra in $M$.
A generic filter may be regarded as a map $G:\Sigma\to 2$.

We use the following facts which are probably all due to Solovay:
\begin{enumerate}
\item \label{product} (see Kunen \cite{K} 3.13)
For any two disjoint sets $\Sigma$ and $\Gamma$ in a countable
standard model $M$,
\begin{enumerate}
 \item $G$ is $\meas(\Sigma\cup \Gamma)$-generic over M iff
 \item $G\res \Sigma$ is
       $\meas(\Sigma)^M$-generic over $M$ and $G\res \Gamma$ is
       $\meas(\Gamma)^{M[G\res \Sigma]}$-generic over $M[G\res \Sigma]$.
\end{enumerate}

\item \label{absolute} (Kunen \cite{K} 3.22)
Suppose $G:\Sigma\to 2$ is $\meas(\Sigma)$-generic over $M$ and
$Y\in M$ is such that
$$M\models Y \su 2^\om \mbox{ has positive outer measure.}$$
Then
$$M[G]\models Y \su \cantor \mbox{ has positive outer measure.}$$

\item \label{lowen} (Well-known)
Suppose $G:\om_2\to 2$ is $\meas(\om_2)$-generic over $M$ and
$$M[G]\models Y\mbox{ has positive outer measure.}$$
Then
there exists a set $\Sigma\su\om_2$ in $M$ of cardinality
$\om_1$ in $M$ such that if $Z=M[G\res \Sigma]\cap Y$, then
$$M[G\res \Sigma]\models Z \mbox{ has positive outer measure.}$$
\end{enumerate}

Fact \ref{lowen} is proved with a Lowenheim-Skolem argument as follows.
Let $f:\cantor\to 2\cross\cantor$ be a
map with the following property:
If $f(x)=(i,z)$, then
\begin{enumerate}
  \begin{enumerate}
    \item $i=1$ iff $x\in Y$ and
    \item if $x$ is a code for a Borel set of measure zero set $Z(x)$,
          then $z\in Y\setminus Z(x)$.
  \end{enumerate}
\end{enumerate}
Since there is a recursive pairing function taking
$2\cross\cantor$ to $\cantor$ it
suffices to show that for any function $f:\cantor\to\cantor$
in $M[G]$ there exists a set $\Sigma\su\om_2$ in $M$
of size $\om_1$ in $M$ such that $\cantor\cap M[G\res \Sigma]$ is
closed under $f$ and $$f\res(M[G\res \Sigma])\in M[G\res \Sigma].$$
For any $x\in\cantor\cap M[G]$ there
exists a sequence $(B_n:n\in\om)$ of Borel sets in $M$ with countable
support such that for any $n\in\om$ we have $x(n)=1$ iff $G\in B_n$
(the equivalence class of $B_n$ is the boolean value of
the statement ``$x(n)=1$'').
Any such sequence $(B_n:n\in\om)$ is called a canonical name for
an element of $\cantor$ (see Kunen \cite{K} 3.17).
Working in the ground model $M$ with a name
for the function $f$, we can define a map $F$ from canonical names
to canonical names such that for any canonical name $\tau$,
$F(\tau)$ will be a canonical name for $f(\tau^G)$.
Since canonical names have countable support and $M$ satisfies
the GCH there exists a set $\Sigma\su\om_2$ of cardinality $\om_1$
in $M$ such that for any canonical name $\tau$ with support from
$\Sigma$, the support of $F(\tau)$ is a subset of $\Sigma$.
This proves Fact \ref{lowen}.

\bigskip
To prove Lemma \ref{randmod}.2, suppose
$Y\su\rr$ has positive outer measure.  By Fact \ref{lowen}
above
there exists a set $\Sigma\su\om_2$ in $M$ of cardinality
$\om_1$ in $M$ such that if $Z=M[G]\cap Y$, then
$$M[G\res \Sigma]\models Z \mbox{ has positive outer measure.}$$
Now since $M$ is a model of CH we have that
$M[G\res \Sigma]$ is a model of CH (see Kunen \cite{K} 3.14).
Hence $Z$ has cardinality $\om_1$.  By Facts
\ref{product} and \ref{absolute}, it follows that
$Z$ has positive outer measure in $M[G]$.
\qed

Finally we prove Theorem \ref{thm3}.  By Lemma \ref{equiv}
if there were such a nonmeasurable function, then
there would be
 a set $H\su \rr^2$ such that
 $\la (H^y )=0$ and $Y=\{ x: \la (\rr \se H_x )=0 \}$ has
 positive outer measure. By applying Lemma \ref{randmod}.2
 we get $Z\su Y$ with positive outer measure
 and cardinality $\om_1$. By Lemma \ref{randmod}.1
 we know that the reals are not covered by the
 $\om_1$ measure zero sets $\{\rr \se H_x:x\in Z\}$.
 Suppose $y\notin \bigcup \{\rr \se H_x:x\in Z\}$.
 Then $y\in\bigcap \{H_x:x\in Z\}$ which implies
 $Z\su H^y$ contradicting the fact that $H^y$ has
 zero measure.
\qed

\bigskip

{\bf Remarks.}
The next statement is implicit in Freiling [8] (see the proof of the Theorem
on p. 198). The following are equivalent:
\begin{description}
\item[{(i)}] there is a function
  $f:[0,1]\times [0,1]\to [0,1]$ such that $f_x , f^y$ are
  measurable for every $x$
  and $y,$ and $\int (\int f_x dy )dx \ne \int (\int f^y dx )dy$;
\item[{(ii)}] there exists a set $H\subset [0,1] \times [0,1]$
  such that $H^y$ is a null set for every $y$ and
  $[0,1]\setminus H_x $ is a null set for every x.
\end{description}
This is similar to our Lemma \ref{equiv};
also, it implies that if Fubini's theorem is not true for arbitrary
bounded functions, then
there is a nonmeasurable function $f$ such that
$f_x$ is approximately continuous and $f^y$ is measurable for every
$x,y.$

\begin{center}
  Addresses
\end{center}

\bigskip
A. Miller:  York University,  Department of Mathematics,
 North York,  Ontario M3J 1P3, Canada (Permanent address:
 University of Wisconsin-Madison,
  Department of Mathematics,
 Van Vleck Hall,
 480 Lincoln Drive,
 Madison, Wisconsin 53706-1388, USA). e-mail: miller@math.wisc.edu

\medskip
M. Laczkovich: E\"{o}tv\"{o}s Lor\'{a}nd University, Department
of Analysis, Budapest, M\'{u}zeum krt. 6-8, H-1088, Hungary.
e-mail: laczk@ludens.elte.hu

\bigskip
November 1994

\end{document}